\documentclass[10pt,twoside]{article}
\usepackage[colorlinks]{hyperref}
\usepackage[margin=1in]{geometry}
\usepackage{times,amsmath,amsthm,amsfonts,amssymb,mathtools,color,latexsym,amsxtra,layout,bibentry,makecell}

\newcommand{\abs}[1]{\left | #1 \right |}
\newcommand{\para}[1]{\left ( #1 \right )}
\newcommand{\bracket}[1]{\left [ #1 \right ]}
\newcommand{\inner}[1]{\left \langle #1  \right \rangle}
\newcommand{\honeinner}[1]{\left \langle #1  \right \rangle_{H^1}}
\newcommand{\htwoinner}[1]{\left \langle #1  \right \rangle_{H^2}}
\newcommand{\hthreeinner}[1]{\left \langle #1  \right \rangle_{H^3}}
\newcommand{\ltwoinner}[1]{\left \langle #1  \right \rangle_{L^2}}
\newcommand{\norm}[1]{\left \| #1 \right \|}
\newcommand{\lpnorm}[1]{\left \| #1 \right \|_{L^p}}

\newcommand{\ltwonorm}[1]{\left \| #1 \right \|_{L^2}}
\newcommand{\lfournorm}[1]{\left \| #1 \right \|_{L^4}}
\newcommand{\linftynorm}[1]{\left \| #1 \right \|_{L^\infty}}
\newcommand{\honenorm}[1]{\left \| #1 \right \|_{H^1}}
\newcommand{\htwonorm}[1]{\left \| #1 \right \|_{H^2}}
\newcommand{\hthreenorm}[1]{\left \| #1 \right \|_{H^3}}
\newcommand{\set}[1]{\left \{ #1 \right \}}

\newcommand{\eqn}[1]{\begin{equation}\begin{split} #1 \end{split}\end{equation}}
\newcommand{\eqnnolabel}[1]{\begin{equation*}\begin{split} #1 \end{split}\end{equation*}}

\newcommand{\R}{\mathbb{R}}

\newcommand{\N}{\mathbb{N}}
\newcommand{\Z}{\mathbb{Z}}

\newcommand{\mx}{\mathcal{X}}

\newcommand{\my}{\mathcal{Y}}
\newcommand{\mh}{\mathcal{H}}

\newcommand{\me}{\mathcal{E}}

\newcommand{\D}[1]{\frac{\partial}{\partial #1}}

\newcommand{\into}{\hookrightarrow}

\newcommand{\too}{\longrightarrow}

\newcommand{\case}[1]{\begin{cases} #1 \end{cases}}
\newcommand{\tr}{\textup{Tr}}

\newcommand{\weaklyto}{\rightharpoonup}
\newcommand{\weaklystarto}{\overset{\ast}{\rightharpoonup}}

\newcommand{\sgn}{\textup{sign}}

\usepackage{color,xfrac}

\theoremstyle{plain}
\newtheorem{theorem}{Theorem}[section]
\newtheorem{lemma}[theorem]{Lemma}
\newtheorem{o-thm}[theorem]{Theorem}
\newtheorem{proposition}[theorem]{Proposition}
\newtheorem{cor}[theorem]{Corollary}

\theoremstyle{definition}

\author{Hagop Karakazian \thanks{Hagop Karakazian, Department of Mathematics, American University of Beirut, Lebanon (hk93@aub.edu.lb)} \and Nabil Nassif \thanks{Nabil Nassif, Department of Mathematics, American University of Beirut, Lebanon (nn12@aub.edu.lb)}}

\title{\bf Strong and Weak Solutions to the Hasegawa-Mima Equation with Periodic Boundary Conditions}
\date{\today}

\begin{document}
\maketitle

\begin{abstract}
\noindent The two dimensional Hasegawa-Mima (HM)
equation
$$ -\Delta u_t+u_t = \{u,\Delta u\} + ku_y$$
describes the time evolution of drift waves in magnetically-confined plasma. Several authors have treated the HM equation theoretically \cite{paumond04,guo04} and numerically \cite{hariri}, with difficulties arising when handling the non-linear Poisson's bracket $\{u,\Delta u\}:=u_x\Delta u_u-u_y\Delta u_x $. In this paper, we introduce a new decoupling approach that avoids the Poisson's bracket term by reformulating the HM equation as a system of two linear PDEs, a solution of which is a pair $(u,w)$ such that
$$\textup{(HM)}\,\,\,\left\{\begin{array}{lll}
w_t + \vec{V}(u) \cdot \nabla w = ku_y\\
-\Delta u+u=w,  \\
\end{array}\right.$$
\noindent where $\vec{V}(u)= -u_y \vec{\textbf{i}} + u_x \vec{\textbf{j}}$ is a divergence-free vector field. Based on this coupled hyperbolic-elliptic system, we derive several variational frames, all propitious  for finding weak solutions with spacial periodic boundary conditions and lower regularity assumptions on the initial data.  More precisely, for initial data $u_0 \in H_P^2(\Omega)$ with $w_0:=(I-\Delta) u_0 \in L^2(\Omega)$, we prove the existence of a weak solution that is global in time. And for initial data $u_0 \in H_P^3(\Omega)$ with $w_0:=(I-\Delta) u_0 \in H_P^1(\Omega) \cap L^\infty(\Omega)$, we prove the existence of a unique strong solution that is local in time. Our proofs are based on the existence of fixed-point ordered pairs  $\set{u_N,w_N}$ that solve Petrov-Galerkin HM systems, constructed using  spacial Fourier basis. Through appropriate a-priori estimates combined with compactness arguments, we reach when $N\to\infty$ limit point solutions $(u,w)$ to the (HM) system.
\end{abstract}

\vspace{0.2in}


\vspace{0.1in}


\section{Introduction}	
Magnetic plasma confinement is one of the most promising ways in future energy production. To understand the phenomena, several mathematical models can be found in literature\cite{iterarticle,hm77,hm78,hw83}, of which the simplest and powerful two dimensional turbulent system model is the Hasegawa-Mima (HM) equation that describes the time evolution of drift waves in magnetically-confined plasma. Originally derived by Akira Hasegawa and Kunioki Mima during late 70s\cite{hm77,hm78}, but when normalized, it can\cite{chm89,hariri} be put as the following PDE that is third order in space and first order in time:
\begin{equation}\label{hm} 
-\Delta u_t+u_t = \{u,\Delta u\} + ku_y
\end{equation}
\noindent where $\{u,v\}=u_xv_y-u_yv_x$ is Poisson's bracket, $u(x,y,t)$ describes the electrostatic potential, $k=\partial_x (\ln \frac{n_0}{\omega_{ci}})$ depends on the background particle density $n_0$ and the ion cyclotron frequency $\omega_{ci}$, which in turn depends on the initial magnetic field. Here, $k=0$ refers to homogeneous plasma. As a cultural note, equation (\ref{hm}) is also referred to as the Charney-Hasegawa-Mima equation in geophysical context that models the time-evolution of Rossby waves in the atmosphere\cite{chm89}.

Using semi-group theory combined with compactness methods, local existence and uniqueness of the Cauchy problem for \eqref{hm} have been obtained in \cite{paumond04} for $u_0 \in H^m(\R^2)$ with $m\geq 4$, and global existence for $u_0 \in H^2(\R^2)$. Similarly, local existence and uniqueness of \eqref{hm} with periodic boundary conditions on a square domain $\Omega$ have been obtained in \cite{karakazian} for $u_0 \in H_P^m(\Omega)$  with $m \geq 4$ (see \eqref{defhmp} for definition of $H_P^m(\Omega)$). 

In this paper, we circumvent the highly non-linear Poisson's bracket $\{u,\Delta u\}$, both theoretically and computationally. We introduce a new variable $w=-\Delta u + u$ and write
$$\{u,\Delta u\} =\{u,u-w\}=\{u,u\}+\{u,-w\}=-\{u,w\}=\{w,u\}=w_xu_y-w_yu_x.$$
This allows us to obtain an equivalent formulation of the HM equation as the following coupled system of hyperbolic-elliptic PDEs:
\begin{equation}\label{hmsystem}
\left\{\begin{array}{lll}
w_t + \vec{V}(u) \cdot \nabla w = ku_y\\
-\Delta u+u=w,  \\
\end{array}\right.
\end{equation}
\noindent where $\vec{V}(u)= -u_y \vec{\textbf{i}} + u_x \vec{\textbf{j}}$ is a divergence-free vector field. This new formulation will be naturally amenable not only to obtain solutions for lower regularity assumptions on the initial data $u_0$, but also to provide a finite element scheme for obtaining a numerical simulation. 

Following \cite{hariri,karakazian}, we also take into account the toroidal shape of plasma confinement chambers, such as Tokamaks, and consider the HM equation on $\Omega=(0,L)\times (0,L)$ with the following periodic boundary conditions: $u \in H^1(\Omega)$ with
$$\mbox{(PBC's)} \case{ u(0,y)=u(L,y) & \mbox{a.e. } y \in [0,L] \\ u(x,0)=u(x,L) & \mbox{a.e. } x \in [0,L]}$$

\noindent where values on the boundary are defined via the trace operator $\tr: H^1(\Omega) \too H^{1/2}(\partial \Omega)$ given by $\tr(u)=u|_{\partial \Omega}$.

\noindent Accordingly, we define a new frame of \textit{periodic Sobolev spaces} as follows:
\eqn{\label{defhmp} H_P^0(\Omega)&=L^2(\Omega) \\
H_P^1(\Omega) &= \set{u \in H^1(\Omega) \ \bigr\vert \  u \mbox{ satisfies PBC's}}, \\
H_P^m(\Omega)&=\set{u \in H^m(\Omega)  \ \bigr\vert \ u,u_{x},u_y \in H_P^{m-1}(\Omega)} \mbox{ for integers } m \geq 2.}
    We remark here that $H_P^m(\Omega)$ is a closed subspace of $H^m(\Omega)$, and so itself is a Hilbert space. A significant consequence of the PBC's is the vanishing integrals over $\partial \Omega$ when integrating by parts, which in turn implies the skew-symmetry of $\partial_x$ and $\partial_y$. More precisely,
$$\int_\Omega f_x g - \int_\Omega fg_x= fg|_{x=L} - fg|_{x=0}=0,  \,\, \forall f,g \in H_P^1(\Omega).$$

In this paper, we consider and solve the following formulations of the HM coupled system \eqref{hmsystem} with PBC's according to different lower regularity assumptions on the initial data:

\subsection*{Strong $H^1$ Formulation of the HM Coupled System}
Given $u_0\in H_P^3(\Omega)$ with $w_0:=u_0-\Delta u_0 \in H_P^1(\Omega) \cap L^\infty(\Omega)$, seek $T>0$ and $(u,w)$ such that:
\begin{equation}\label{hmstrong}
\left\{\begin{array}{lll}
w_t + \vec{V}(u) \cdot \nabla w = ku_y & \mbox{in } L^\infty(0,T;L^2(\Omega)) &(.1)\\
-\Delta u+u=w &\mbox{in }  H_P^1(\Omega), \forall t \in (0,T)  &(.2)\\
\mbox{PBC's on } u,\, u_x,\, u_y,\,w & \mbox{ a.e. on }  \partial \Omega\times (0,T) & (.3)\\
u(0)=u_0 \mbox{ and }w(0)=w_0 & \mbox{ a.e. on } \overline{\Omega}  & (.4)\\
\end{array}\right.
\end{equation} 
Taking the $L^2$ inner-product of equations in (\ref{hmstrong}.1) and (\ref{hmstrong}.2) with a test function $v\in H^1_P(\Omega)$ leads us to an equivalent H\textsuperscript{1} variational form:
\begin{equation}\label{hmh1}
\left\{\begin{array}{ll}
\ltwoinner{w_t,v} + \ltwoinner{\vec{V}(u) \cdot \nabla w,v}= \ltwoinner{ku_y,v}, \, \forall v \in H^1_P(\Omega), \ \forall t \in (0,T) &(.1)\\
\honeinner{u,v}=\ltwoinner{w,v},\forall v \in H^1_P(\Omega), \ \forall t \in (0,T) &(.2)\\
w(0)=w_0 & (.3)\\
u(0)=u_0 &(.4)\\
\end{array}\right.
\end{equation}

\subsection*{Weak $L^2$ Formulation of the HM Coupled System}

\noindent Due to the divergence-free property of $\vec{V}(u)$, Green's formula yields
\eqn{\label{eqnskew2} \ltwoinner{\vec{V}(u) \cdot \nabla w,v}= - \ltwoinner{\vec{V}(u) \cdot \nabla v, w},\,\forall\,v,w\,\in H_P^1(\Omega),}
which for $v=w$, gives us an essential identity 
\eqn{\label{eqnskew} \ltwoinner{\vec{V}(u) \cdot \nabla w,w}=0, \, \forall w \in H_P^1(\Omega).}
Identity (\ref{eqnskew2}) and integration of (\ref{hmh1}.1) over the temporal interval $[0,t]$, with $0<t<T$, allow us to relax spacial and temporal regularity requirements on $(u,w)$ and reach to the following weak $L^2$ formulation of the HM coupled system: \\

\noindent Given $u_0\in H_P^2(\Omega)$ with $w_0:=u_0-\Delta u_0 \in L^2(\Omega)$, seek $(u,w)$ such that:
\begin{equation}\label{hml2}
\left\{\begin{array}{ll}
\ltwoinner{w(t)-w_0,v} = \displaystyle \int_0^t \ltwoinner{\vec{V}(u(s)) \cdot \nabla v,w(s)} + \ltwoinner{k u_y(s),v} ds,\forall v \in H^2_P(\Omega), \ \forall t \in (0,T] &(.1)\\
\honeinner{u,v}=\ltwoinner{w,v},\forall v \in H^1_P(\Omega), \ \forall t \in (0,T] &(.2)\\
w(0)=w_0 & (.3)\\
u(0)=u_0 & (.4) \\
\end{array}\right.
\end{equation}
Formulation \eqref{hml2} has been used in
\cite{nsk} for obtaining a robust algorithm for simulating the HM model.

\subsection*{Main Results} \label{summary}

Our main results are the following two theorems:

\begin{theorem} \label{thmexist5}
Given initial data $u_0 \in H_P^2(\Omega)$ with $w_0:=(I-\Delta) u_0$, the weak $L^2$ formulation of the HM coupled system \eqref{hml2} has a solution pair $(u,w)$ with the following regularities:
$$u \in L^\infty(0,T; H_P^2(\Omega))\cap C([0,T],H_P^1(\Omega))$$
and
$$w \in L^\infty(0,T; L^2(\Omega)),$$
where the existence time 
\eqn{\label{Tind} T:= (C_E\linftynorm{k}+1)^{-1},}
\noindent independent of the initial data, allowing us to extend $(u,w)$ globally to $T=\infty$. Here, $C_E>0$ is the elliptic regularity constant from the elliptical equation. 
\end{theorem}

\noindent Theorem \ref{thmexist5} is proved in section \ref{secproof1}. 


\begin{theorem} \label{thmexist4}
Given initial data $u_0 \in H_P^3(\Omega)$ with $w_0:=(I-\Delta) u_0 \in H_P^1(\Omega) \cap L^\infty(\Omega)$, the strong $H^1$ formulation of the HM coupled system \eqref{hmstrong} has a unique local solution pair $(u,w)$ on $[0,T]$, where the existence time $T=T(u_0,w_0)$ is inversely dependent on $(u_0,w_0)$ as in equation \eqref{T}, with the  following regularities for the solution:
$$u \in L^\infty(0,T; H_P^3(\Omega)) \cap C([0,T], H_P^2(\Omega)) \mbox{ with } u_t \in L^\infty(0,T; H_P^2(\Omega)),$$
and
$$w \in L^\infty(0,T; H_P^1(\Omega)\cap L^\infty(\Omega)) \cap C([0,T], L^2(\Omega)) \mbox{ with }  w_t \in L^\infty(0,T; L^2(\Omega)).$$
\end{theorem}

\noindent Theorem \ref{thmexist4}  is proved in section \ref{secproof2}. \\

\noindent As a byproduct and solely for the hyperbolic equation:
\begin{enumerate}
    \item [(i)] Global existence is obtained in Corollary \ref{Exist-Hyper-L2}  for $w_0 \in L^2(\Omega)$, and,
    \item [(ii)] Global existence and uniqueness are obtained in Corollary \ref{Exist-Hyper-H1} for $w_0 \in H^1_P(\Omega) \cap L^\infty(\Omega)$.
\end{enumerate}

\section{Methodology}

To establish our main results, we study the elliptic and hyperbolic equations, each on its own. Our results are summarized as follows:
\begin{proposition} \label{Reg-Elliptic}
If $w \in H_P^m(\Omega)$, with $m\geq 0$, then the elliptic equation (\ref{hmstrong}.2) has a unique strong solution $u=\mathcal{E}(w):=(I-\Delta)^{-1}w \in H_P^{m+2}$ in the sense that:
\begin{equation}\label{Res-E}
\left\{\begin{array}{ll}
-\Delta u + u = w & \mbox{ in } H_P^m(\Omega) \\
\norm{u}_{H^{m+2}} \leq C_E \norm{w}_{H^m} &
\end{array}\right.
\end{equation}
\noindent where $C_E>0$ is an elliptic regularity constant depending only on $m$. We note here that the temporal regularity of $w$ carries over to $u$.
\end{proposition}

\noindent This result is well-known for Dirichlet boundary value problems. However for the sake of completeness, we provide a proof in Appendix \ref{Sec-Elliptic} using an equivalent Hilbert-Fourier basis $\set{\phi_j/\sqrt{\lambda_j^m}}_{j=1}^\infty$ for $H_P^m(\Omega)$, where $\phi_j$'s and $\lambda_j$'s are the eigenvectors and and eigenvalues (in non-decreasing order) of $I-\Delta$ on $L^2(\Omega)$ satisfying PBC's on $\partial\Omega$.

We note here that the operator $(I-\Delta)$ is also invertible over all finite dimensional subspaces $$E_N:=\mbox{span}\set{\phi_1, \phi_2, \hdots, \phi_N}$$ of $H_P^1(\Omega)$. Namely if $w=\sum_{j=1}^N c_j\phi_j \in E_N$, then $u=\sum_{j=1}^N \lambda_j^{-1}c_j\phi_j \in E_N$. Motivated by this fact, we give a Galerkin formulation of the hyperbolic equation (\ref{hmh1}.1) on $E_N$,
and establish its existence and uniqueness in the following proposition.

\begin{proposition} \label{thmexist1}
Given $u\in C([0,T],H_P^2(\Omega))$, then for every $N \in \mathbb{Z}^+$, there exists a unique function\linebreak $w_N:=\mh_N(u) \in C^1((0,T),E_N) \cap C([0,T],E_N)$ satisfying:
\begin{equation}\label{Disc-Hyper}
\left\{\begin{array}{ll}
\ltwoinner{w_N',v} = \ltwoinner{\vec{V}(u) \cdot \nabla v,w_N} + \ltwoinner{ku_y,v}  \, \forall v \in E_N, \, \forall t \in (0,T) & (.1) \\ 
w_N(0)=(I-\Delta)\mbox{proj}_{E_N}u(0) & (.2)
\end{array}\right.
\end{equation} 
\noindent where $'$ denotes the derivative with respect to $t$.

\end{proposition}

\begin{proof}
Write $w_N=\sum_{i=1}^N c_i(t)\phi_i$ and set $v=\phi_j$ in (\ref{Disc-Hyper}.1) for $j=1,\cdots,N$, then
$$c_j'(t) + \sum_{i=1}^N \underbrace{\ltwoinner{\vec{V}(u(t)) \cdot \nabla \phi_i,\phi_j}}_{A_{i,j}(t)}c_i(t) =\underbrace{\ltwoinner{ku_y(t),\phi_j}}_{\vec{F}_j(t)} \, \forall j=1,\cdots,N.$$
Setting $\vec{C}_{0,j}:=c_j(0)=\ltwoinner{\mbox{proj}_{E_N}u_0,\phi_j}$, this becomes the system of $N$ ODE's:
$$\vec{C}'(t) + A(t)\vec{C}(t)=\vec{F}(t) \mbox{ \ with \ } \vec{C}(0)=\vec{C}_0.$$
Since for $u\in C([0,T],H_P^2(\Omega))$,  $\{A_{ij}(t)\}$ and $\{F_j(t)\}$ are defined and continuous on $[0,T]$, then by a usual Picard iteration, we get a unique solution $\vec{C}(t) \in C^1((0,T),\R^N) \cap C([0,T],\R^N)$, which completes the proof.

\end{proof}
\noindent We will call $\mh_N$ the Galerkin solution operator for the hyperbolic equation ($w_N=\mh_N(u)$), and consider a sequence of fixed-point problems, each consisting of finding a fixed-point $u_N$ for $\me\circ\mh_N$ ($u_N=\me\circ\mh_N(u_N)$), visualized by:

\begin{table}[h]
\centering
\begin{tabular}{|l|l|l|}
\cline{1-1} \cline{3-3}
 & $\xrightarrow[\textup{of the Galerkin hyperbolic PDE}]{\textup{Solver  operator } \mh_N \textup{ }}$ &  \\
$u_N \in \mx_N$ &  & $w_N \in \my_N$  \\
 & $\xleftarrow[\textup{of the elliptic PDE}]{\textup{\quad \quad 
 Solver Operator }\me \quad \quad}$ &  \\ \cline{1-1} \cline{3-3} 
\end{tabular}
\end{table}

\noindent where $\mx_N$ and $\my_N$ are non-empty, closed, bounded, convex subsets of $C([0,T], E_N)$ that will be specified according to the regularity assumptions on the initial data $u_0$ and $w_0$. Once specified, we will have $\me$ continuous and compact, and each $\mh_N$ continuous, thus making each $\me \circ \mh_N$ a continuous compact mapping, which according to Schauder's theorem will have a fixed-point $u_N=\me \circ \mh_N(u_N)$ satisfying \eqref{Res-E} and \eqref{Disc-Hyper}.\\

\noindent Via this setup, we derive a priori estimates for the sequence of pairs $(u_N,w_N)$ in selected Bochner Hilbert spaces. Due to uniform boundedness and compact embeddings,  we then extract a subsequence of fixed-point pairs and show they converge in the weak-star topology in time and weak topology in space to the solution pairs $(u,w)$ claimed in Theorems \ref{thmexist5} and \ref{thmexist4}.


\section{Proof of Theorem \ref{thmexist5}} \label{secproof1}

\subsection*{A Priori $L^2$ Estimate for the Galerkin Solution}

\begin{lemma} \label{lemma1}
The Galerkin solution $w_N=\mathcal{H}_N(u)$ associated to $u\in C([0,T],H_P^2(\Omega))$ satisfies the following estimate:
\eqn{\label{estimate1}
\norm{w_N}_{L^\infty(L^2)} \leq \linftynorm{k} \norm{u}_{L^\infty(H^1)} T+ \ltwonorm{w_0}
}
\end{lemma}

\begin{proof}
Setting $v=w_N(t)\in E_N$ in (\ref{Disc-Hyper}.1) and using (\ref{eqnskew}) we obtain
$$\frac{1}{2}\frac{d}{dt}\ltwonorm{w_N}^2=\int_\Omega k u_y w_N,$$
\noindent which via Cauchy-Schwarz inequality becomes
$$\ltwonorm{w_N}\frac{d}{dt}\ltwonorm{w_N} \leq \linftynorm{k} \norm{u}_{L^\infty(H^1)} \ltwonorm{w_N}$$
Cancelling out $\ltwonorm{w_N}$, integrating over $[0,t]$, with $0<t\leq T$, and taking the sup over $[0,T]$ completes the proof.
\end{proof}
\subsection*{Setup of the Sequence of Fixed-Point Problems}

Given an initial data $u_0 \in H_P^2(\Omega)$ with $w_0:=(I-\Delta)u_0$, we pick $T:= (C_E\linftynorm{k}+1)^{-1}$ as in \eqref{Tind}, and set
\eqn{\label{eqntl2} C_\mx=\frac{C_E \ltwonorm{w_0}}{1-C_E\linftynorm{k}T}\geq0,}
\noindent and for each $N \in \Z^+$, consider $\mx_N$ and $\my_N$ to be the following non-empty, closed, bounded, convex sets:
\eqnnolabel{\mx_N:=\set{u \in C([0,T], E_N) \quad | \quad \norm{u}_{L^\infty(H^2)} \leq C_\mx} \\ \my_N:=\set{w \in C([0,T], E_N) \quad | \quad C_E \norm{w}_{L^\infty(L^2)} \leq C_\mx}}

\noindent Observe that each solution operator is well-defined due to Propositions \ref{Reg-Elliptic} and \ref{thmexist1}, and the inclusions:
\begin{enumerate}
\item [(i)] $\mh_N(\mx) \subseteq \my$ through estimate (\ref{estimate1}) as
\eqnnolabel{C_E \norm{\mh_N(u)}_{L^\infty(L^2)} &\leq C_E\linftynorm{k}T\norm{u}_{L^\infty(H^1)} + C_E\ltwonorm{w_0} \\
& \leq C_E\linftynorm{k}TC_\mx + C_E \ltwonorm{w_0} =C_\mx.}
\item [(ii)] $\me (\my) \subseteq \mx$ through estimate (\ref{Res-E}) as
$$\norm{\me (w)}_{L^\infty(H^2)} \leq C_E \norm{w}_{L^\infty(L^2)} \leq C_\mx.$$
\end{enumerate}

\noindent Also due to finite-dimensionality of $E_N$ that rises the norm-equivalence
$$\ltwonorm{v}  \leq \honenorm{v} \leq \sqrt{\lambda_N} \ltwonorm{v} \quad \forall v \in E_N$$
on $E_N$, we obtain the continuity of $\mh_N$ as follows.

\begin{lemma} \label{contofhnl2}
$\mh_N:\mx_N \too \my_N$ is continuous.
\end{lemma}

\begin{proof}
Let $\epsilon >0$ be given, set $\delta=\epsilon/(C_L^2\lambda_NC_\mx T / C_E +\linftynorm{k}T)$, and
let $u_1, u_2 \in \mx_N$ be such that $\norm{u_1-u_2}_{L^\infty(H^2)} < \delta$. Subtracting (\ref{Disc-Hyper}.1) with $w_2=\mh_N(u_2)$ from that of $w_1=\mh_N(u_1)$ and setting $v=w_1-w_2$, we obtain, via (\ref{eqnskew}),
$$\frac{1}{2}\frac{d}{dt}\ltwonorm{w_1-w_2}^2 = \ltwoinner{\vec{V}(u_2-u_1) \cdot \nabla w_1, w_2-w_1} + \ltwoinner{k(u_1-u_2)_y,w_2-w_1}.$$
Now applying H\"older's inequality to the left hand side and cancelling $\ltwonorm{w_1-w_2}$ from both sides, we obtain

\eqnnolabel{\frac{d}{dt}\ltwonorm{w_1-w_2} &\leq \lfournorm{\vec{V}(u_2-u_1)} \lfournorm{\nabla w_1} +\linftynorm{k} \ltwonorm{(u_1-u_2)_y} \\
&\leq C_L^2\htwonorm{u_1-u_2} \htwonorm{w_1} +\linftynorm{k} \honenorm{u_1-u_2} \\
&\leq (C_L^2\lambda_N\norm{w_1}_{L^2} +\linftynorm{k}) \norm{u_1-u_2}_{H^2}\\
&< (C_L^2\lambda_N C_\mx/C_E +\linftynorm{k})\delta = \epsilon/T.}
\noindent Integrating this over $[0,t]$, with $t \leq T$, and taking the sup over $[0,T]$, we finally obtain
$$\norm{\mh_N(u_1) - \mh_N(u_2)}_{L^\infty(L^2)} = \norm{w_1-w_2}_{L^\infty(L^2)} < \epsilon.$$
\end{proof}

\subsection*{Obtaining a Candidate Solution Pair $(u,w)$}

At this point, due to Schauder fixed-point theorem, we have a sequence of fixed-point pairs $\set{u_N,w_N}$ that are uniformly bounded in $L^\infty(0,T,H_P^2(\Omega)) \times L^\infty(0,T,L^2(\Omega))$ satisfying \eqref{Res-E} and \eqref{Disc-Hyper}. Moreover:
\begin{lemma}
The sequence of temporal derivative $u_N'$ is also uniformly bounded in $L^\infty(0,T;L^2(\Omega))$ as
\eqn{\label{estimate5}
\norm{u_N'}_{L^\infty(L^2)} \leq C_E[4C_L^2 (\linftynorm{k} C_\mx T+ \ltwonorm{w_0})+\linftynorm{k}]C_\mx.}


\end{lemma}
\begin{proof}
For convenience, we drop the index $N$ from $u_N$ and $w_N$. Setting $v=(I-\Delta)^{-1}u'$ in (\ref{Disc-Hyper}.1) and using the self-adjointness of $(I-\Delta)$ on $L^2(\Omega)$, we obtain
\eqnnolabel{\ltwonorm{u'}^2 &= \ltwoinner{\vec{V}(u) \cdot \nabla [(I-\Delta)^{-1}u'], w}+\ltwoinner{ku_y,(I-\Delta)^{-1}u'} \\
&\leq \norm{\vec{V}(u) }_{L^4} \norm{\nabla [(I-\Delta)^{-1}u']}_{L^4} \ltwonorm{w}  + \linftynorm{k} \ltwonorm{u_y} \ltwonorm{\nabla [(I-\Delta)^{-1}u']} \\
&\leq [4C_L^2 \ltwonorm{w} +\linftynorm{k}]\htwonorm{u} \htwonorm{(I-\Delta)^{-1}u'} \\
&\leq C_E[4C_L^2 \ltwonorm{w} +\linftynorm{k}]\htwonorm{u} \ltwonorm{u'},
}

\noindent Finally, simplifying and taking the sup over $[0,T]$ and using estimate \eqref{estimate1}, estimate \eqref{estimate5} follows.
\end{proof}

\noindent Thus the sequence  $\set{u_N}$ is uniformly bounded in
$$\set{v \in L^\infty(0,T;H_P^2(\Omega)) \ | \ v_t \in L^\infty(0,T;L^2(\Omega))} \subset \subset  C([0,T],H_P^1(\Omega)).$$
This last compact embedding is due to Aubin-Lions lemma, which in itself is based on Rellich's Theorem \cite{rellich1}.

Now along with Banach–Alaoglu theorem, we can eventually extract a subsequence of pairs $(u_N, w_N)$ (re-indexed for convenience) such that: 

\eqn{\label{uNweakstarh2} u_N \weaklystarto \mbox{ some } u \in L^\infty(0,T;H_P^2(\Omega)), \mbox{ i.e } \int_0^T \htwoinner{u_N,v} dt \too \int_0^T \htwoinner{u,v} dt  \quad \forall v \in L^1(0,T;H_P^2(\Omega)),}

\eqn{\label{wNweakstarl2} w_N \weaklystarto \mbox{ some }  w \in L^\infty(0,T;L^2(\Omega)), \mbox{ i.e } \int_0^T \ltwoinner{w_N,v} dt \too \int_0^T \ltwoinner{w,v} dt  \quad \forall v \in L^1(0,T;L^2(\Omega)),}

\eqn{\label{uNstrongl2} u_N \too u \in C([0,T],H_P^1(\Omega)), \mbox{ i.e } \norm{u_N-u}_{L^
\infty(H^1)} \too 0.}

\noindent Note that $u$ in equation \eqref{uNstrongl2} is equal to that of \eqref{uNweakstarh2} due to uniqueness of weak star limits in $L^\infty(0,T;H^1(\Omega))$.

\subsection*{Passing to the Limit}

We now show that this pair $(u,w)$ is a solution to the weak $L^2$ form of the HM coupled system \eqref{hml2}. 

We have $-\Delta u_N+u_N=w_N$ where the RHS converges to $w$, and the LHS converges to $-\Delta u + u$, both weakly star in $L^\infty(0,T; L^2(\Omega))$. Thus uniqueness of weak star limits implies that $-\Delta u + u=w$ in $L^\infty(0,T; L^2(\Omega))$.

Regarding the hyperbolic equation, let $t \in (0,T]$ and $v \in H_P^2(\Omega)$. Then there exist a sequence $\{v_M\}$ with $v_M \in E_M$ such that $\norm{v_M}_{H^2} \leq \norm{v}_{H^2}$ and $v_M \to v$ strongly in $H_P^2(\Omega)$. Now fix $M\geq 1$, then for every $N \geq M$, integrating (\ref{Disc-Hyper}.1) over $[0,t]$, we get (due to $E_M \subset E_N$)
\eqn{\label{h1hyperaux}
\ltwoinner{w_N(t)-w_N(0),v_M}=  \int_0^t  \ltwoinner{\vec{V}(u_N(s)) \cdot \nabla v_M,w_N(s)} ds + \int_0^t \ltwoinner{ku_{N,y}(s),v_M} ds}
\noindent where as $N \to \infty$, we have:
\begin{enumerate}
\item [(i)] $w_N(0)$ converges strongly to $w_0$ in $L^2(\Omega)$ by definition.

\item [(ii)] $w_N$ converges to $w$ weakly in $L^2(0,t; L^2(\Omega)) \ \forall t \in (0,T]$ via \eqref{wNweakstarl2} with $v=\chi_{[0,t]\times \Omega} \widetilde{v}$ where $\widetilde{v} \in L^2(0,T;L^2(\Omega))$.

\item [(iii)] A subsequence (depending on $v_M$) of the continuous functions $f_N(t):=\ltwoinner{w_N(t),v_M}$ converges uniformly to $\ltwoinner{w(t),v_M}$ on $[0,T]$ by Arzela-Ascoli theorem. To see this, observe that $f_N(t)$ is uniformly bounded on $[0,T]$ as
$$\abs{f_N(t)} \leq C_\mx \norm{v_M}_{L^2}/C_E,$$
and equicontinuous on $(0,T]$ as
\eqnnolabel{\abs{f_N(t)-f_N(s)} &=\abs{\ltwoinner{\int_s^t w_N'(\tau) \ d\tau, v_M}} \\
&\overset{\mbox{Fubini's}}{=}\abs{\int_s^t \ltwoinner{w_N',v_M} \ d \tau } \\
&\leq \abs{\int_s^t \ltwoinner{\vec{V}(u_N) \cdot \nabla v_M, w_N} \ d \tau } + \abs{\int_s^t \ltwoinner{ku_{N,y}, v_M} \ d \tau }  \\
&\leq (t-s) \bracket{\norm{\vec{V}(u_N)}_{L^\infty(L^4)} \norm{\nabla v_M}_{L^4} \norm{w_N}_{L^\infty(L^2)} + \linftynorm{k}\norm{u_{N,y}}_{L^\infty(L^2)} \norm{v_M}_{L^2}} \\
&\leq (t-s) \bracket{C_L^2\norm{u_N}_{L^\infty(H^2)}\norm{v_M}_{H^2} \norm{w_N}_{L^\infty(L^2)} + \linftynorm{k}\norm{u_N}_{L^\infty(H^1)}\norm{v_M}_{L^2}} \\
&\leq (t-s) \bracket{C_L^2C_\mx^2/C_E + \linftynorm{k}C_\mx}\norm{v}_{H^2}}
\noindent where $C_L>0$ is the constant in the Ladyzhenskaya's inequality
$$\lfournorm{v} \leq C_L \ltwonorm{v}^{1/2} \cdot \honenorm{v}^{1/2}  \quad \forall v \in H^1(\Omega).$$
\noindent Hence a subsequence $\{f_{N_k}(t)\}$ converges uniformly to $\ltwoinner{\widetilde{w}(t),v_M}$ on $(0,T]$. By the uniqueness of weak limits in $L^2(0,t;L^2(\Omega))$, 
$\widetilde{w}=w$.

\item [(iv)] Re-indexing according to the subsequence found in part (iii), we have $$\label{carefulconv1} \int_0^t \ltwoinner{\vec{V}(u_N(s)) \cdot \nabla v_M, w_N(s)} ds \too \int_0^t \ltwoinner{\vec{V}(u(s)) \cdot \nabla v_M, w(s)} ds$$
as
\eqnnolabel{& \abs{\int_0^t \ltwoinner{\vec{V}(u_N(s)) \cdot \nabla v_M, w_N(s)} ds - \int_0^t \ltwoinner{\vec{V}(u(s)) \cdot \nabla v_M, w(s)} ds} \\
&\hspace{0.2in} \leq \underbrace{\abs{\int_0^t \ltwoinner{\vec{V}(u_N(s)-u(s)) \cdot \nabla v_M, w_N(s)} ds}}_{A_1} + \underbrace{\abs{\int_0^t \ltwoinner{\vec{V}(u(s)) \cdot \nabla v_M, w_N(s)-w(s)} ds}}_{B_1},}
\noindent where
$$A_1 \leq \underbrace{\norm{\nabla v_M}_{L^\infty(L^\infty)}}_{\textup{fixed}}  \underbrace{\norm{w_N}_{L^2(L^2)}}_{\leq \textup{constant}}  \underbrace{\norm{\vec{V}(u_N(s)-u(s))}_{L^2(L^2)}}_{\too 0 \textup{ by } \eqref{uNstrongl2}} \too 0,$$
and $B_1\too 0$ by the weak convergence $w_N \weaklyto w$ in $L^2(0,t; L^2(\Omega))$.

\item [(v)] Finally through $u_N \weaklyto u$ in $L^2(0,t;H_P^1(\Omega))$ or strongly in $C([0,T];H_P^1(\Omega))$, we have
$$\int_0^t \ltwoinner{ku_{N,y}(s),v_M} ds \to \int_0^t \ltwoinner{ku_y(s), v_M} ds.$$
\end{enumerate}

\noindent Putting all together, for each $M \geq 1$, equation \eqref{h1hyperaux} converges to
\eqn{\label{convergenceM} \ltwoinner{w(t)-w_0,v_M}=  \int_0^t  \ltwoinner{\vec{V}(u(s)) \cdot \nabla v_M ,w(s)} ds + \int_0^t \ltwoinner{ku_y(s),v_M} ds \quad \mbox{ for all } t \in (0,T],}
where the choice of the subsequence depends on $M$, but ultimately, equation \eqref{convergenceM} holds for all $M$. Finally letting $M \to \infty$, \eqref{convergenceM} converges to (\ref{hml2}.1).

\subsection*{Initial Conditions}

Finally, we show that the initial conditions are satisfied. Taking the limit of (\ref{hml2}.1) as $t \to 0$, we obtain
$$\lim_{t \to 0} \ltwoinner{w(t)-w_0,v}=0 \quad \forall v \in H_P^2(\Omega),$$
which, via the elliptical equation and the self-adjointness of $I-\Delta$, is equivalent to 
$$\lim_{t \to 0} \ltwoinner{u(t)-u_0,v}=0  \quad \forall v \in L^2(\Omega).$$

This completes the proof of Theorem \ref{thmexist5}. As a byproduct, we obtain the existence sole for the hyperbolic equation:

\begin{cor} \label{Exist-Hyper-L2} For every $T>0$, if $u \in C([0,T],H_P^2(\Omega))$ and $w_0 \in L^2(\Omega)$, then the hyperbolic equation has a weak solution $w\in L^\infty(0,T;L^2(\Omega))$ in the sense of (\ref{hml2}.1) and (\ref{hml2}.3).
\end{cor}

\begin{proof}
Mimic proof of Theorem \ref{thmexist5} as follows: Directly after a priori estimate \eqref{estimate1}, extract the subsequence $w_N$ and go to ``Passing to the Limit'' section with $u_N=u$.
\end{proof}

\section{Proof of Theorem \ref{thmexist4}}
\label{secproof2}

\subsection*{A Priori $H^1$ Estimates for the Galerkin Solution}

\begin{lemma} \label{lemma2} If $u \in C([0,T],H_P^2(\Omega)) \cap L^\infty(0,T; H_P^3(\Omega))$ with $w_0 \in H_P^1(\Omega) \cap L^\infty(\Omega)$, then
\eqn{\label{estimate2} \norm{w_N}_{L^\infty(L^\infty)} &\leq 21C_\infty \norm{k}_{W^{2,\infty}} \norm{u}_{L^\infty(H^3)}T+ 2\linftynorm{w_0},}
\eqn{\label{estimate3} \norm{w_N}_{L^\infty(H^1)} \leq 168C_\infty\norm{k}_{W^{2,\infty}}\bracket{T \norm{u}_{L^\infty(H^3)}}^2 + (16\linftynorm{w_0} + 2\linftynorm{k})\bracket{T \norm{u}_{L^\infty(H^3)}}+ \honenorm{w_0},}
\noindent and 
\eqn{\label{estimate4} \norm{w_N'}_{L^\infty(L^2)} \leq (2 C_\infty \norm{w_N}_{L^\infty(H^1)} +\linftynorm{k}) \norm{u}_{L^\infty(H^3)},}
\noindent where $C_\infty>0$ is the constant from the embedding $H^2(\Omega) \into L^\infty(\Omega)$.
\end{lemma}

\begin{proof} For convenience, we will write $w$ instead of $w_N$. \\

\noindent \textbf{Derivation of \eqref{estimate2}}:
Let $p\geq 4$ be an even integer. Setting $v=\mbox{proj}_{E_N}[pw^{p-1}(t)]$ in (\ref{Disc-Hyper}.1) and using Lemma \ref{propsetup} (i) on the left-hand-side term, we obtain 
\eqn{\label{aux1} \frac{d}{dt} \lpnorm{w}^p=\underbrace{- \ltwoinner{ \vec{V}(u) \cdot \nabla w,  \mbox{proj}_{E_N}[pw^{p-1}]}}_{Z} + \ltwoinner{\mbox{proj}_{E_N}[ku_y],pw^{p-1}},} 
\noindent where
$$Z= - \ltwoinner{u_x, w_y\mbox{proj}_{E_N}[pw^{p-1}]} + \ltwoinner{u_y, w_x\mbox{proj}_{E_N}[pw^{p-1}]}.$$
Now applying Lemma \ref{propsetup} (iv) on $w^{p-1}$, \, $w_y\mbox{proj}_{E_N}[w^{p-1}]$ and $w_x\mbox{proj}_{E_N}[w^{p-1}]$, there exists a large enough integer $Q \geq N$ such that 
$$w_y\mbox{proj}_{E_N}[w^{p-1}]= \mbox{proj}_{E_Q}[w_yw^{p-1}] \textup{ \ and \ } w_x\mbox{proj}_{E_N}[w^{p-1}]= \mbox{proj}_{E_Q}[w_xw^{p-1}],$$
so that
\eqnnolabel{Z &=- \ltwoinner{\mbox{proj}_{E_Q}[u_x], \D{y}w^p}+\ltwoinner{\mbox{proj}_{E_Q}[u_y],\D{x}w^p} \\
&=\ltwoinner{\D{y}\mbox{proj}_{E_Q}[u_x]-\D{x}\mbox{proj}_{E_Q}[u_y],w^p}\\
&=\ltwoinner{\mbox{proj}_{E_Q}[u_{xy}]-\mbox{proj}_{E_Q}[u_{yx}],w^p}=0.}

\noindent Hence (\ref{aux1}) implies
$$ \lpnorm{w}^{p-1} \frac{d}{dt}\lpnorm{w} \leq  \linftynorm{\mbox{proj}_{E_N}[ku_y]} \norm{w}_{L^{p-1}}^{p-1},$$
where
$$\linftynorm{\mbox{proj}_{E_N}[ku_y]} \leq C_\infty \norm{ku_y}_{L^\infty(H^2)} \leq 21 C_\infty\norm{k}_{W^{2,\infty}} \norm{u}_{L^\infty(H^3)}$$
and
$$\norm{w}_{L^{p-1}}^{p-1} \leq L^{2/p} \lpnorm{w}^{p-1}$$
using Lemma \ref{propsetup} (ii) and 
H\"older's inequality, respectively. Cancelling out $\lpnorm{w}^{p-1}$, we obtain
$$\frac{d}{dt}\lpnorm{w} \leq 21 C_\infty \norm{k}_{W^{2,\infty}}\norm{u}_{L^\infty(H^3)} L^{2/p}.$$
Integrating this over $[0,t]$ gives
$$\lpnorm{w(t)} \leq 21 C_\infty \norm{k}_{W^{2,\infty}}   \norm{u}_{L^\infty(H^3)} L^{2/p} t + L^{2/p}\linftynorm{w(0)}.$$

\noindent Finally, letting $p \to \infty$ and then taking the sup over $[0,T]$, estimate \eqref{estimate2} follows. \\

\noindent \textbf{Derivation of \eqref{estimate3}}:
Setting $v=(I-\Delta)w$ in (\ref{Disc-Hyper}.1) and using (\ref{eqnskew}), we obtain
\eqn{\label{aux3} \frac{1}{2}\frac{d}{dt} \honenorm{w}^2 &= -\ltwoinner{\nabla [\vec{V}(u) \cdot \nabla w], \nabla w} + \ltwoinner{k\nabla u_y, \nabla w},}
\noindent where $-\ltwoinner{\nabla [\vec{V}(u) \cdot \nabla w], \nabla w}$ simplifies through integration by parts to 
\eqnnolabel{&\ltwoinner{u_{xxx}w_y,w} +\ltwoinner{u_{xx}w_{yx},w} +\ltwoinner{u_{xyy}w_{y},w} +\ltwoinner{u_{xy}w_{yy},w} \\
&-\ltwoinner{u_{xyx}w_x,w} - \ltwoinner{u_{xy}w_{xx},w} -\ltwoinner{u_{yyy}w_x,w} -\ltwoinner{u_{yy}w_{xy},w}}
\noindent whose terms can be bounded above by $\linftynorm{w}\hthreenorm{u} \honenorm{w}$ directly or through Lemma \ref{propsetup} (v). Thus after cancelling $\honenorm{w}$ out, (\ref{aux3}) implies
$$\frac{d}{dt} \honenorm{w} \leq (8 \norm{w}_{L^\infty(L^\infty)}+ 2\linftynorm{k})\norm{u}_{L^\infty(H^3)}.$$
Now integrating over $[0,t]$, with $t \leq T$, and taking the sup over $[0,T]$, we obtain
$$\norm{w}_{L^\infty(H^1)} \leq  (8 \norm{w}_{L^\infty(L^\infty)}+ 2\linftynorm{k})\norm{u}_{L^\infty(H^3)} T + \honenorm{w(0)}.$$
Finally via estimate (\ref{estimate2}) and $\honenorm{w(0)} \leq \honenorm{w_0}$, estimate \eqref{estimate3} follows. \\

\noindent \textbf{Derivation of \eqref{estimate4}}:
Setting  $v=w'$ in (\ref{Disc-Hyper}.1), we obtain
\eqnnolabel{\ltwonorm{w'}^2&=-\ltwoinner{u_yw_x,w'}+\ltwoinner{u_xw_y,w'}+\ltwoinner{ku_y, w'} \\
&\leq  (\linftynorm{u_y} \ltwonorm{w_x}  + \linftynorm{u_x} \ltwonorm{w_y} +\linftynorm{k}\cdot\ltwonorm{u_y}) \ltwonorm{w'},}
\noindent so that
$$\ltwonorm{w'} \leq (2C_\infty \honenorm{w} +\linftynorm{k})\hthreenorm{u}$$
Finally taking the sup over $[0,T]$, estimate \eqref{estimate4} follows.
\end{proof}




\subsection*{Setup of the Sequence of Fixed-Point Problems}

Given a non-zero initial data $u_0 \in H_P^3(\Omega)$ with $w_0:=(I-\Delta) u_0 \in H_P^1(\Omega) \cap L^\infty(\Omega)$, we consider:
\eqnnolabel{\mx_N:=\set{u \in C([0,T], E_N) \quad | \quad T\norm{u}_{L^\infty(H^3)} \leq 1} \\ \my_N:=\set{w \in C([0,T], E_N) \quad | \quad C_ET\norm{w}_{L^\infty(H^1)} \leq 1}}

\noindent where $T$ is the least existence time given by
\eqn{\label{T} T:=\bracket{C_E (168C_\infty\norm{k}_{W^{2,\infty}} + 16\linftynorm{w_0} + 2\linftynorm{k} + \honenorm{w_0}+1)}^{-1}.}

Observe that each solution operator is well-defined because:
\begin{enumerate}
\item [(i)] $\mh_N(\mx_N) \subseteq \my_N$ through estimate (\ref{estimate3}) as
$$C_E T\norm{\mh_N(u)}_{L^\infty(H^1)} \leq C_ET(168C_\infty\norm{k}_{W^{2,\infty}} + 16\linftynorm{w_0} + 2\linftynorm{k}+ \honenorm{w_0}) \leq 1.$$

\item [(ii)]$\me (\my_N) \subseteq \mx_N$ through estimate (\ref{Res-E}) as
$$T\norm{\me (w)}_{L^\infty(H^3)} \leq C_E T \norm{w}_{L^\infty(H^1)} \leq 1.$$
\end{enumerate}

\noindent Moreover,

\begin{lemma} \label{contofhnh1}
$\mh_N:\mx_N \too \my_N$ is continuous.
\end{lemma}

\begin{proof}
Let $\epsilon >0$ be given, set $\delta=\epsilon/[C_L^2\sqrt{\lambda_N} /C_E +\linftynorm{k}T]$, and
let $u_1, u_2 \in \mx_N$ be such that $\norm{u_1-u_2}_{L^\infty(H^3)} < \delta$. Continuing similarly as in the proof of Lemma \ref{contofhnl2}, we arrive to
\eqnnolabel{\frac{d}{dt}\ltwonorm{w_1-w_2} &\leq C_L^2\htwonorm{u_1-u_2} \htwonorm{w_1} +\linftynorm{k} \honenorm{u_1-u_2} \\
&\leq (C_L^2\sqrt{\lambda_N}\norm{w_1}_{H^1} +\linftynorm{k}) \norm{u_1-u_2}_{H^3}\\
&< [C_L^2\sqrt{\lambda_N} /(C_ET) +\linftynorm{k}]\delta = \epsilon/T.}
\noindent Integrating this over $[0,t]$, with $t \leq T$, and taking the sup over $[0,T]$, we finally obtain
$$\norm{\mh_N(u_1) - \mh_N(u_2)}_{L^\infty(H^1)} = \norm{w_1-w_2}_{L^\infty(H^1)} \leq \sqrt{\lambda_N} \norm{w_1-w_2}_{L^\infty(L^2)} < \epsilon.$$
\end{proof}

\subsection*{Obtaining a Candidate Solution Pair $(u,w)$}

At this point, due to Schauder fixed-point theorem, we have a sequence of fixed-point pairs $\set{u_N,w_N}$ that are uniformly bounded in
$$L^\infty(0,T;H_P^3(\Omega)) \times L^\infty(0,T;H_P^1(\Omega) \cap L^\infty(\Omega)),$$
with their temporal derivatives $\set{u_N',w_N'}$ uniformly bounded in 
$$L^\infty(0,T;H_P^2(\Omega)) \times L^\infty(0,T;L^2(\Omega)),$$
satisfying \eqref{Res-E} and \eqref{Disc-Hyper}. And so by Banach–Alaoglu theorem, we can eventually extract a subsequence of pairs (re-indexed for convenience) such that:

\eqn{\label{uNweakstarh3} u_N \weaklystarto \mbox{ some } u \in L^\infty(0,T;H_P^3(\Omega)), \mbox{ i.e } \int_0^T \hthreeinner{u_N,v} dt \too \int_0^T \hthreeinner{u,v} dt  \quad \forall v \in L^1(0,T;H_P^3(\Omega)),}

\eqn{\label{wNweakstarh1} w_N \weaklystarto \mbox{ some }  w \in L^\infty(0,T;H^1(\Omega)), \mbox{ i.e } \int_0^T \honeinner{w_N,v} dt \too \int_0^T \honeinner{w,v} dt  \quad \forall v \in L^1(0,T;H^1(\Omega)),}

\eqn{\label{wNprimeweakstarlinfty} w_N \weaklystarto w \in L^\infty(0,T;L^\infty(\Omega)), \mbox{ i.e } \int_0^T \int_\Omega w_Nv \too \int_0^T \int_\Omega wv  \quad \forall v \in L^1(0,T;L^1(\Omega)),}

\eqn{\label{uNprimeweakstarh1} u_N' \weaklystarto u_t \in L^\infty(0,T;H^2(\Omega)),\mbox{ i.e } \int_0^T \htwoinner{u'_N,v} dt \too \int_0^T \htwoinner{u_t,v} dt \quad \forall v \in L^1(0,T;H^2(\Omega)).}

\eqn{\label{wNprimeweakstarl2} w_N' \weaklystarto w_t \in L^\infty(0,T;L^2(\Omega)),\mbox{ i.e } \int_0^T \ltwoinner{w'_N,v} dt \too \int_0^T \ltwoinner{w_t,v} dt \quad \forall v \in L^1(0,T;L^2(\Omega)),}

\noindent Note that $w$ in  \eqref{wNprimeweakstarlinfty} is  equal to that of \eqref{wNweakstarh1} due to uniqueness of weak limits in $L^2(0,T;L^2(\Omega))$. Also $u_t$ and $w_t$ in  \eqref{uNprimeweakstarh1} and \eqref{wNprimeweakstarl2} are the temporal derivatives of $u$ and $w$ in \eqref{uNweakstarh3} and \eqref{wNweakstarh1}, respectively, due to uniqueness of weak   derivatives. Also due to Aubin-Lions lemma, we have the compact embeddings:
$$\set{v \in L^\infty(0,T;H_P^3(\Omega)) \ | \ v_t \in L^\infty(0,T;H_P^2(\Omega)) } \subset \subset C([0,T],H_P^2(\Omega)),$$
$$\set{v \in L^\infty(0,T;H_P^1(\Omega)) \ | \ v_t \in L^\infty(0,T;L^2(\Omega)) } \subset \subset C([0,T];L^2(\Omega)),$$
through which we obtain the following strong convergences:
$$u_N \too u \in C([0,T];H_P^2(\Omega)), \mbox{ i.e } \norm{u_N-u}_{L^\infty(H^2)}\too 0$$
$$w_N \too w \in C([0,T];L^2(\Omega)), \mbox{ i.e } \norm{w_N-w}_{L^
\infty(L^2)}\too 0.$$

\subsection*{Passing to the Limit}

We now show that this pair $(u,w)$ is a solution to the Hasegawa-Mima Coupled System (\ref{hmstrong}).

We have $-\Delta u_N+u_N=w_N$ where the RHS converges to $w$, and the LHS converges to $-\Delta u + u$, both weakly star in $L^\infty(0,T; H_P^1(\Omega))$. Thus uniqueness of weak star limits implies that $-\Delta u + u=w$ in $L^\infty(0,T; H_P^1(\Omega))$.

Regarding the hyperbolic equation, we let $t \in (0,T]$ and $v \in L^2(\Omega)$ and consider a sequence $\{v_M\}$ with $v_M \in E_M$ such that $\ltwonorm{v_M} \leq \ltwonorm{v}$ and $v_M \to v$ strongly in $L^2(\Omega)$. Now fix $M \geq 1$, then for every $N \geq M$, integrating (\ref{Disc-Hyper}.1) over $[0,t]$, we get (due to $E_M \subset E_N$ and \eqref{eqnskew2}),
\eqn{\label{h1hyperaux2}
\int_0^t \ltwoinner{w_N',v_M} ds+ \ltwoinner{\vec{V}(u_N) \cdot \nabla w_N,v_M} ds = \int_0^t \ltwoinner{ku_{N,y},v_M} ds}
\noindent Now using the weak convergence of $w_N'$ to $w_t$ in $L^2(0,t;L^2(\Omega))$, and the claims (i), (ii), (iv) and (v) in the ``Passing to the Limit'' section in the proof of Theorem \ref{thmexist5}, equation \eqref{h1hyperaux2} converges to
$$\int_0^t \ltwoinner{w_t,v_M} + \ltwoinner{\vec{V}(u) \cdot \nabla w,v_M} ds= \int_0^t \ltwoinner{ku_y,v_M} ds$$
as $N \to \infty$, with $w(0)=w_0$. Letting $M \to \infty$, we obtain
\eqn{\label{final1} \int_0^t \ltwoinner{w_t+\vec{V}(u)\cdot \nabla w- ku_y,v} ds =0 \, \, \, \forall v \in L^2(\Omega) \mbox{ and } \forall t \in (0,T]}
Subtracting equation \eqref{final1} with $t-h$ instead of $t$ from the original, dividing both sides by $h$, and letting $h \to 0$, Lesbegue's Differentiation Theorem implies
\eqn{\label{final2}  \ltwoinner{w_t(t)+\vec{V}(u(t))\cdot \nabla w(t)- ku_y(t),v} =0 \, \, \, \forall v \in L^2(\Omega) \mbox{ and a.e. } t \in(0,T]}
\noindent Finally, setting $v=w_t(t)+\vec{V}(u(t)) \cdot \nabla w(t) - ku_y(t) \in L^2(\Omega)$ in \eqref{final2}, equation (\ref{hmstrong}.1).

\subsection*{Uniqueness of the Solution Pair $(u,w)$}

Finally suppose  $(u_1,w_1)$ and $(u_2,w_2)$ are two solution pairs in the context of Theorem \ref{thmexist4}. Then subtracting the hyperbolic equation with $(u_2,w_2)$ from that of $(u_1,w_1)$, and taking the $L^2$-innerproduct of the resulting equation with $w_1-w_2$, we obtain
$$\frac{1}{2}\frac{d}{dt}\ltwonorm{w_1-w_2}^2  = \underbrace{\ltwoinner{\vec{V}(u_1-u_2) \cdot \nabla (w_1-w_2), w_1}}_{A_2}+ \underbrace{\ltwoinner{k(u_1-u_2)_y,w_1-w_2}}_{B_2},$$
where by Cauchy-Schwarz inequality and Lemma \ref{propsetup} (v),
\eqnnolabel{A_2 &\leq \linftynorm{w_1} \para{\ltwoinner{\abs{(u_1-u_2)_x},\abs{(w_1-w_2)_y}} + \ltwoinner{\abs{(u_1-u_2)_y},\abs{(w_1-w_2)_x}}} \\
&\leq \linftynorm{w_1} \para{\ltwoinner{\abs{(u_1-u_2)_{xy}},\abs{w_1-w_2}} + \ltwoinner{\abs{(u_1-u_2)_{yx}},\abs{w_1-w_2}}} \\
&\leq 2\linftynorm{w_1} \cdot \htwonorm{u_1-u_2} \cdot \ltwonorm{w_1-w_2} \\
&\leq 2C_E\linftynorm{w_1} \cdot \ltwonorm{w_1-w_2}^2,}
\noindent and
$$B_2 \leq \linftynorm{k}\cdot\htwonorm{u_1-u_2} \cdot \ltwonorm{w_1-w_2}  \leq \linftynorm{k}C_E \ltwonorm{w_1-w_2}^2.$$

\noindent Now putting all together and cancelling the term $\ltwonorm{w_1-w_2}$, we obtain
\eqnnolabel{\frac{d}{dt}\ltwonorm{w_1-w_2}&\leq \underbrace{(2\norm{w_1}_{L^\infty(L^\infty)}+ \linftynorm{k})  C_E}_{\leq \textup{constant}}\ltwonorm{w_1-w_2},}
\noindent where $\ltwonorm{w_1-w_2}$ is continuous in $t$ and $w_1(0)=w_2(0)$, so that by Gronwall's inequality, we have $\ltwonorm{w_1(t)-w_2(t)}=0$ where $t$ was arbitrary. Thus, $w_1=w_2$ in $C([0,T], L^2(\Omega))$, and via the uniqueness of the solution to the elliptic equation, $u_1=u_2$ in $C([0,T], H_P^2(\Omega))$.

This completes the proof of Theorem \ref{thmexist4}. As a byproduct, we obtain the existence sole for the hyperbolic equation:

\begin{cor} \label{Exist-Hyper-H1} For every $T>0$, if $u \in C([0,T],H_P^2(\Omega)) \cap L^\infty(0,T; H_P^3(\Omega))$ and $w_0 \in H^1_P(\Omega) \cap L^\infty(\Omega)$, then the hyperbolic equation has a unique strong solution 
$$w \in L^\infty(0,T; H_P^1(\Omega)\cap L^\infty(\Omega)) \cap C([0,T], L^2(\Omega)) \mbox{ with }  w_t \in L^\infty(0,T; L^2(\Omega)),$$
satisfying (\ref{hmstrong}.1) and $w(0)=w_0$.
\end{cor}

\begin{proof}
Mimic proof of Theorem \ref{thmexist4} as follows: Directly after a priori estimates \eqref{estimate2}, \eqref{estimate3} and \eqref{estimate4}, extract the subsequences involving $w_N$ and $w_N'$, then go to ``Passing to the Limit'' section with $u_N=u$.
\end{proof}

\section{Concluding Remarks}

Thus far, we have shown global weak existence for the $L^2$ formulation \eqref{hml2} to and unique local existence for the $H^1$ formulation \eqref{hmstrong} to the HM equation using a priori estimates and compactness methods to extract convergent subsequences to the sequence of fixed-points $u_N=\me \circ \mh_N(u_n)$. To our knowledge, the following problems remain open:
\begin{enumerate}
    \item Uniqueness of a global weak solution pair $(u,w)$ in Theorem \ref{thmexist5}.
    \item Existence of a global strong solution pair $(u,w)$ in Theorem \ref{thmexist4}.
    \item Uniqueness of a global weak solution $w$ solely for the hyperbolic equation in Corollary \ref{Exist-Hyper-L2}.
\end{enumerate}

\noindent On the computational side, the procedures used in this paper can be implemented using spectral methods.

\newpage

\appendix
\section{Appendix}

\subsection{Proof of Proposition \ref{Reg-Elliptic} Using an Equivalent Hilbert-Fourier Basis for $H_P^m(\Omega)$}\label{Sec-Elliptic}

We begin with a lemma that allows us to obtain a Hilbert-Fourier basis for $L^2(\Omega)$ that satisfy PBC's. For that purpose, let $H_P^{-1}(\Omega)$ be the dual of $H_P^1(\Omega)$ and denote its action on $H_P^1(\Omega)$ by $\inner{\cdot,\cdot}_{H^{-1},H^1}$. Then,

\begin{lemma}
The operator $(I-\Delta)^{-1}:H_P^{-1}(\Omega) \too H_P^1(\Omega)$ exists and is self-adjoint and compact on $L^2(\Omega)$.
\end{lemma}

\begin{proof}
By Riesz-Fr\'echet representation theorem, for every $w \in H_P^{-1}(\Omega)$, there exists a unique $u\in H_P^1(\Omega)$, such that
\eqn{\label{defnhdualinner} \inner{w,v}_{H^{-1},H^1}=\honeinner{u,v}=\ltwoinner{u,v} + \ltwoinner{\nabla u , \nabla v} \quad \forall v \in H_P^1(\Omega)}
\noindent with
\eqn{\label{defnhdualnorm} \norm{w}_{H^{-1}}=\honenorm{u}.}

\noindent Observe that \eqref{defnhdualinner} is the variational formulation of $(I-\Delta)u=w \in H_P^{-1}(\Omega)$ through Green's formula and use of PBC's on $u$. And so, the operator $(I-\Delta): H_P^1(\Omega) \too H_P^{-1}(\Omega)$ is invertible with $u=(I-\Delta)^{-1}w$.

Now using \eqref{defnhdualnorm} and Cauchy-Schwarz inequality, it can be shown directly that $\norm{w}_{H^{-1}} \leq \ltwonorm{w}$ whenever $w \in L^2(\Omega)$. And so we have
$$(I-\Delta)^{-1}|_{L^2(\Omega)}: L^2(\Omega) \xhookrightarrow[]{\textup{continuously}} H_P^{-1}(\Omega) \xrightarrow[]{(I-\Delta)^{-1}} H_P^1(\Omega)  \xhookrightarrow[]{\textup{compactly}} L^2(\Omega)$$
where the compact embedding is due to Rellich's theorem \cite{rellich1}. Hence $(I-\Delta)^{-1}$ is compact on $L^2(\Omega)$.

Considering the context of the triplet $H_P^1(\Omega) \subset L^2(\Omega) \simeq L^2(\Omega)^* \subset H_P^{-1}(\Omega)$, where all inclusions are continuous and dense, we have that
\eqn{\inner{w,v}_{H^{-1},H^1}=\ltwoinner{w,v} \quad \forall w \in L^2(\Omega), \ \forall v \in H_P^1(\Omega).}
\noindent This allows us to show that $(I-\Delta)^{-1}$ is self-adjoint on $L^2(\Omega)$ as
\eqn{\label{selfdual}\ltwoinner{w,(I-\Delta)^{-1}v} = \inner{w,(I-\Delta)^{-1}v}_{H^{-1},H^1}= \honeinner{(I-\Delta)^{-1}w,(I-\Delta)^{-1}v} \quad \forall w,v \in L^2(\Omega)}
and the last term is symmetric. 
\end{proof}

Now by Theorem 6.11 of \cite{brezis}, $L^2(\Omega)$ has a Hilbert basis consisting of eigenvectors $\set{\phi_j}_{j=1}^\infty$ of $(I-\Delta)^{-1}$ with eigenvalues $\set{\eta_j}_{j=1}^\infty$. Now since the eigenvalues of an invertible operator cannot be zero, then
\eqn{(I-\Delta)^{-1}\phi_j=\eta_j \phi_j \iff (I-\Delta)\phi_j=\frac{1}{\eta_j} \phi_j }
\noindent so that $\set{\phi_j}_{j=1}^\infty$ are eigenvectors of $I-\Delta$ with eigenvalues $\lambda_j:=1/\eta_j$.
It's not hard to see that \eqn{\phi_j=\frac{e^{2\pi i \vec{\bold{x}} \cdot \vec{\xi}_j /L}}{L} \ \textup{ and } \lambda_j=1+\frac{4\pi^2\abs{\vec{\xi}_j }^2}{L^2} \mbox{ with }  \vec{\xi}_j \in \N \times \N,}

\noindent and all $\phi_j$'s satisfy the PBC's. This basis allows us to characterize the spaces $H_P^m(\Omega)$ as follows.

\begin{proposition} \label{propelliptic}
Let $m \geq 0$ be an integer. There exists an innerproduct $\inner{\inner{\cdot,\cdot}}_m$ on $H_P^m(\Omega)$ with

\eqn{\label{propdeltaij1} \inner{\inner{\phi_i,\phi_j}}_m = \lambda_i^m \delta_{ij} = \case{\lambda_i^m & \mbox{for }i=j \\ 0 & \mbox{for }i\not=j}.}
\noindent Thus making $\set{\phi_i/\sqrt{\lambda_i^m}}_{i=1}^\infty$ a Hilbert basis for $H_P^m(\Omega)$ and inducing the Hilbert norm

\eqn{\label{propdeltaij2} \abs{\norm{\cdot}}^2_m=\sum_{i=1}^\infty \lambda_i^m \ltwoinner{\cdot, \phi_i}^2}
that is equivalent to the usual Sobolev norm:
\eqn{\label{propequivenorm} \norm{\cdot}_{H^m} \leq \abs{\norm{\cdot}}_m \leq (1+\kappa_m) \norm{\cdot}_{H^m}}
where $\kappa_m>0$ is a constant depending on $m$.
\end{proposition}

\begin{proof}
Suppose $m \geq 1$ and consider the sequence $$Q_m: H_P^{2m-2}(\Omega) \too L^2(\Omega)$$ of differential operators given recursively as follows
\eqn{\case{Q_1=0 & \\ Q_m=Q_{m-1}(I-\Delta) - \sum_{\abs{\alpha}\leq m-2} (-1)^{\abs{\alpha}} D^\alpha D^\alpha \Delta & m \geq 2,}}
\noindent where $\alpha=(\alpha_1, \alpha_2)$ is a bi-index of length $\abs{\alpha}=\alpha_1+ \alpha_2$ and $D^\alpha = \partial_x^{\alpha_1} \partial_y^{\alpha_2}$. 

It is worth to note here that the domain of $Q_m$ can extended by density to all of $L^2(\Omega)$ as $H_P^{2m-2}(\Omega)$ contains a Fourier basis of $L^2(\Omega)$.

\begin{lemma}
$Q_m$ is a positive self-adjoint operator on $L^2(\Omega)$, and so it has a positive self-adjoint square root 
 $$Q_m^{\sfrac{1}{2}}: H_P^{m-1}(\Omega) \too L^2(\Omega)$$ 
such that $\para{Q_m^{\sfrac{1}{2}}}^2=Q_m$.
\end{lemma}

\begin{proof}
Self-adjointness of $Q_m$ follows immediately from $\partial^2_{xy}=\partial^2_{yx}$ and the self-adjointness of $\partial^2_{xx}$, $\partial^2_{xy}$ and $\partial^2_{yy}$. Now to show that $Q_m$ is positive, we proceed by induction as follows: For $f \in H_P^{2m-2}(\Omega)$,

Base cases: $\ltwoinner{Q_1 f, f}=0 \geq 0$ and $\ltwoinner{Q_2 f, f}= - \ltwoinner{\Delta f, f}=\ltwonorm{\nabla f}^2 \geq 0$

 Inductive step: Assume that $Q_{m-1}$ is positive for some $m\geq 2$, then
$$\ltwoinner{Q_m f, f} =  \ltwoinner{Q_{m-1} f, f} + \ltwoinner{Q_{m-1} \nabla f, \nabla f} + \norm{\nabla f}_{H^{m-2}}^2 \geq 0.$$
Finally, as $Q_m$ is a differential operator that lowers spacial regularity by $2(m-1)$, then $Q_m^{\sfrac{1}{2}}$ will be a differential operator that lowers spacial regularity by $m-1$. For example: $Q_2^{\sfrac{1}{2}}=\nabla$ and $Q_3^{\sfrac{1}{2}}=2\nabla(I-\Delta)^{\sfrac{1}{2}}$.
\end{proof}

\noindent \textit{Proof of Proposition \ref{propelliptic} continued}.
Consider the innerproducts:
\eqn{\case{\inner{\inner{\cdot,\cdot}}_0 := \ltwoinner{\cdot,\cdot} & \mbox{ on } L^2(\Omega) \\
\inner{\inner{\cdot,\cdot}}_1 := \honeinner{\cdot,\cdot} & \mbox{ on } H_P^1(\Omega) \\
\inner{\inner{\cdot,\cdot}}_m := \inner{\cdot,\cdot}_{H^m} + \ltwoinner{Q_m^{\sfrac{1}{2}} \ \cdot, Q_m^{\sfrac{1}{2}} \ \cdot} & \mbox{ on } H_P^m(\Omega) \mbox{ for }  m \geq 2}}
\noindent and denote the induced norm by
$$\abs{\norm{\cdot}}_m := \sqrt{\inner{\inner{\cdot,\cdot}}_m}.$$
Since $Q_m^{\sfrac{1}{2}}$ lowers spacial regularity by $m-1$, then there exists $\kappa_m >0$ such that
$$\ltwonorm{Q_m^{\sfrac{1}{2}} \ \cdot} \leq \kappa_m \norm{\cdot}_{H^{m-1}},$$
and so (\ref{propequivenorm}) follows.

Observe that $\inner{\inner{\phi_i,\phi_j}}_0 = \ltwoinner{\phi_i,\phi_j}=\lambda_i^0 \delta_{ij}$ and 
$$\inner{\inner{\phi_i,\phi_j}}_1=\ltwoinner{\phi_i,\phi_j}+\ltwoinner{\nabla \phi_i,\nabla \phi_j} =\ltwoinner{(I-\Delta)\phi_i, \phi_j} = \lambda_i \ltwoinner{\phi_i,\phi_j} = \lambda_i^1 \delta_{ij}.$$
Now for $m\geq 2$, we have
\eqnnolabel{\inner{\inner{\phi_i,\phi_j}}_m&= \inner{\phi_i,\phi_j}_{H^m} + \ltwoinner{Q_m \phi_i, \phi_j} \\
&=\inner{\phi_i,\phi_j}_{H^m} + \inner{\nabla \phi_i, \nabla \phi_j}_{H^{m-2}} + \lambda_i\ltwoinner{Q_{m-1}\phi_i,\phi_j} \\
&=\inner{\phi_i,\phi_j}_{H^{m-1}} + \inner{\nabla \phi_i, \nabla \phi_j}_{H^{m-1}}+ \lambda_i\ltwoinner{Q_{m-1}\phi_i,\phi_j} \\
&=\inner{\phi_i,\phi_j}_{H^{m-1}} + (\lambda_i-1) \inner{ \phi_i, \phi_j}_{H^{m-1}}+ \lambda_i\ltwoinner{Q_{m-1}\phi_i,\phi_j} \\
&=\lambda_i \inner{\inner{\phi_i,\phi_j}}_{m-1}
}
\noindent from which by induction (\ref{propdeltaij1}) follows. Moreover, because $\set{\phi_i/\sqrt{\lambda_i^m}}_{i=1}^\infty$ spans $L^2(\Omega)$ which contains a dense subset of $H_P^m(\Omega)$, then $\set{\phi_i/\sqrt{\lambda_i^m}}_{i=1}^\infty$ is a Hilbert basis for $H_P^m(\Omega)$.

Now for $f \in H_P^m(\Omega)$, write
$$\sum_{i=1}^\infty \frac{1}{\lambda_i^m} \inner{\inner{f,\phi_i}}_m \phi_i \overset{H_P^m}{=\joinrel=} f \overset{L^2}{=\joinrel=} \sum_{i=1}^\infty \ltwoinner{f,\phi_i} \phi_i,$$
so that by orthogonality of $\phi_i$'s, we have
$$\inner{\inner{f,\phi_i}}_m = \lambda_i^m \ltwoinner{f,\phi_i}$$
from which (\ref{propdeltaij2}) follows by Parseval's identity.

\end{proof}

Finally, to finish the proof of Proposition \ref{Reg-Elliptic}, we write for $w\in H_P^m(\Omega)$, with $m \geq 0$, 
$$w\overset{L^2}{=\joinrel=}\sum_{i=1}^\infty \ltwoinner{w,\phi_i}\phi_i = \sum _{i=1}^\infty (\ltwoinner{u,\phi_i}+\ltwoinner{\nabla u,\nabla \phi_i})\phi_i = \sum_{i=1}^\infty \ltwoinner{u,(I-\Delta)\phi_i}\phi_i = \sum_{i=1}^\infty \lambda_i \ltwoinner{u,\phi_i}\phi_i,$$
so that by orthogonality of $\phi_i$'s, we have
\eqn{\label{components} \lambda_i \ltwoinner{u,\phi_i}=\ltwoinner{w,\phi_i} \quad \forall i\geq 1.}
Hence
$$\abs{\norm{u}}_{m+2}^2 = \sum_{i=1}^\infty \lambda_i^{m+2}\ltwoinner{u,\phi_i}^2 = \sum_{i=1}^\infty \lambda_i^m \ltwoinner{w,\phi_i}^2 = \abs{\norm{w}}_m^2,$$
so that $$u \in H_P^{m+2}(\Omega) \mbox{ \ and \ } \norm{u}_{H^{m+2}} \leq (1+\kappa_m)\norm{w}_{H^m}.$$

\noindent Setting $C_E:=(1+\kappa_m)$ completes the proof.

\subsection{Auxiliary Results For the a Priori Estimates  \eqref{estimate5}, \eqref{estimate2},  \eqref{estimate3}, and \eqref{estimate4}} \label{appendix2}
\begin{lemma} \label{propsetup}
The following assertions hold:
\begin{enumerate}
\item [(i)] For all $u,v \in L^2(\Omega)$, we have $\ltwoinner{\mbox{proj}_{E_N}u,v} =\ltwoinner{u,\mbox{proj}_{E_N}v}$.
\item [(ii)] For all $u \in H_P^m(\Omega)$, with $m \geq 0$, we have $\norm{\mbox{proj}_{E_N} u}_{H^m} \leq \norm{u}_{H^m}$.
\item [(iii)] For all $u \in H_P^2(\Omega)$. we have $\mbox{proj}_{E_N}[(I-\Delta)u]=(I-\Delta)\mbox{proj}_{E_N} u$.
\item [(iv)] For all $N,M \in \Z^+$, there exists an integer $Q \geq M,N$ such that $$uv \in E_Q \mbox{ for all } u \in E_M \mbox{ and } v \in E_N.$$
\item [(v)] For all $u \in H_P^1(\Omega)$ and $\varphi \in E_N$, we have
\eqn{\label{bypartsabseqn}\int_\Omega \abs{u} \cdot \abs{D\varphi}  \leq \int_\Omega \abs{Du} \cdot \abs{\varphi},}
\noindent where $D$ denotes $\partial_x$ or $\partial_y$. By density, the results also holds for $\varphi \in H_P^1(\Omega)$.
\end{enumerate}
\end{lemma}

\begin{proof}
Parts (i)-(iii) follow from direct computations. Part (iv) is based on the fact that the product of two Fourier basis functions can be written as the sum of multiples of Fourier basis functions with larger eigenvalues. For part (v), observe that $u$ and $Du$ are in $L^1_{\textup{loc}}(\Omega)$, so that $\abs{u} \in L^1_{\textup{loc}}(\Omega)$ and is weakly differentiable with $D\abs{u} = \sgn(u)Du$. When $\varphi \in E_N$, we compute
$$\int_\Omega \abs{u} \cdot \abs{D\varphi} = \int_{\set{D\varphi \geq 0}} \abs{u} D\varphi - \int_{\set{D\varphi < 0}} \abs{u} D\varphi.$$
Since $D\varphi$ is a finite sum of Fourier basis functions, and so the domains $\set{D\varphi \geq 0}$ and $\set{D\varphi < 0}$ have finite number of connected components with piecewise smooth boundaries satsifying
\eqn{\label{Dvarphiboundary} \partial \set{D\varphi \geq 0} \setminus \partial \set{D\varphi < 0}  \subset \partial \Omega.}
Now integrating by parts, we obtain
\eqnnolabel{\int_\Omega \abs{u} \cdot \abs{D\varphi} &= \int_{\partial\set{D\varphi \geq 0}} \abs{u} D\varphi  -  \int_{\set{D\varphi \geq 0}} (D \abs{u}) \varphi -\int_{\partial\set{D\varphi < 0}} \abs{u} D\varphi  +  \int_{\set{D\varphi < 0}} (D \abs{u}) \varphi \\
&\leq \int_{\partial\set{D\varphi \geq 0}} \abs{u} D\varphi  -\int_{\partial\set{D\varphi < 0}} \abs{u} D\varphi  +  \int_{\set{D\varphi \geq 0}} \abs{D \abs{u}} \cdot \abs{\varphi} +  \int_{\set{D\varphi < 0}} \abs{D \abs{u}} \cdot \abs{\varphi} \\
&\overset{(\ref{Dvarphiboundary})}{\leq} \int_{\partial\Omega} \abs{u} D\varphi  + \int_\Omega \abs{D \abs{u}} \cdot \abs{\varphi},}
\noindent where the boundary integral vanishes as $\abs{u}$ and $D\varphi$ satisfy PBC's, and so \eqref{bypartsabseqn} follows.

When $\varphi \in H_P^1(\Omega)$, there is a sequence $\set{\varphi_N}$ such that $\varphi_N \in E_N$ and $\varphi_N \too \varphi$ in $H_P^1(\Omega)$. Now we compute 
\eqnnolabel{\int_\Omega \abs{u} \cdot \abs{D \varphi} & \leq \int_\Omega \abs{u} \cdot \abs{D \varphi_N} +  \int_\Omega \abs{u} \cdot \abs{D (\varphi-\varphi_N)} \\
&\overset{\eqref{bypartsabseqn}}{\leq}  \int_\Omega \abs{Du} \cdot \abs{\varphi_N - \varphi + \varphi} +  \int_\Omega \abs{u} \cdot \abs{D (\varphi-\varphi_N)} \\
&\leq \int_\Omega \abs{Du} \cdot \abs{\varphi} + \honenorm{u} \ltwonorm{\varphi_N-\varphi} +  \ltwonorm{u} \honenorm{\varphi-\varphi_N}.}
Letting $N \too \infty$, the result follows.
\end{proof}

\newpage
\bibliographystyle{ieeetr}
\bibliography{references}

\end{document}